\documentclass[reqno]{amsart}

\usepackage{color}

\newcommand{\pp}{ {\partial} }

\newcommand{\A}{\alpha }

\newcommand{\OO}{{\mathcal O}}

\newcommand{\JJ}{{\mathcal J}}

\newcommand{\R} {\mathbb R}
\newcommand{\cuad}{{\sqcap\kern-.68em\sqcup}}

\newcommand{\foral}{\quad\mbox{for all}\quad}

\newcommand{\be}{\begin{equation}}
\newcommand{\ee}{\end{equation}}

\newtheorem{theorem}{Theorem}
\newtheorem{prop}{Proposition}[section]
\newtheorem{corollary}{Corollary}[section]
\newtheorem{remark}{Remark}[section]
\newcommand{\bremark}{\begin{remark} \em}
\newcommand{\eremark}{\end{remark} }

\numberwithin{equation}{section}

\def\cb{\color{blue}}
\def\cb{}

\title{Nonlocal Minimal Lawson Cones}

\author{Juan D\'avila}
\address{\noindent J. D\'avila -
Departamento de Ingenier\'{\i}a Matem\'atica and CMM, Universidad
de Chile, Casilla 170 Correo 3, Santiago, Chile.}
\email{jdavila@dim.uchile.cl}

\author{Manuel del Pino}
\address{\noindent M. del Pino- Departamento de
Ingenier\'{\i}a  Matem\'atica and CMM, Universidad de Chile, Casilla
170 Correo 3, Santiago, Chile.} \email{delpino@dim.uchile.cl}

\author{Juncheng Wei}
\address{\noindent J. Wei -
Department Of Mathematics, Chinese University Of Hong Kong,
Shatin, Hong Kong, and Department of Mathematics, University of British Columbia, Vancouver, B.C., Canada, V6T 1Z2.} \email{wei@math.cuhk.edu.hk}

\begin{document}

\begin{abstract}
We prove the existence of the analog of Lawson's minimal cones for a notion of nonlocal minimal surface introduced by Caffarelli, Roquejoffre and Savin, and establish their stability/instability in low dimensions. In particular we find that there are nonlocal stable minimal cones in dimension 7, in contrast with the case of classical minimal surfaces.
\end{abstract}

\maketitle

\section{Introduction}
In \cite{caffarelli-roquejoffre-savin}, Caffarelli,  Roquejoffre and Savin introduced a nonlocal notion of perimeter
of a set $E $, which generalizes the  $(N-1)$-dimensional surface area of $\pp E$.  For $0<s<1$,
the $s$-perimeter of $E \subset\R^N$ is defined {\cb (formally)} as
$$
Per_s(E) =  \int_E \int_{\R^N \setminus E}  \frac {dx\, dy}{ |x-y|^{N+ s} }.
$$
This notion is localized to a bounded open set $\Omega$ by setting
 $$Per_s(E,\Omega) = \int_E \int_{\R^N \setminus E}\frac {dx\, dy}{ |x-y|^{N+ s} }  - \int_{E\setminus \Omega} \int_{\R^N \setminus (E\cup\Omega)}  \frac {dx\, dy}{ |x-y|^{N+ s} }. $$  This quantity makes sense, even if the last two terms above are infinite, by rewriting it in the form
$$
Per_s(E,\Omega) =  \int_{E\cap\Omega} \int_{\R^N \setminus E} \frac {dx\, dy}{ |x-y|^{N+ s} }  +  \int_{E\setminus\Omega} \int_{\Omega \setminus E} \frac {dx\, dy}{ |x-y|^{N+ s} }  .
$$
 Let us assume that $E$ is an open set  set with $\pp E\cap \Omega$ smooth.  The usual notion of perimeter is recovered by the formula
 \begin{equation}
\label{1new}
 \lim_{s\to 1}  (1-s) Per_s(E,\Omega) =   Per(E,\Omega) = c_N {\mathcal H}^{N-1} (\pp E \cap \Omega),
 \end{equation}
 see \cite{savin-valdinoci}. Let us consider a unit normal vector field $\nu$  of $\Sigma = \pp E$ pointing to the exterior of $E$,
  and consider functions $h\in C_0^\infty (\Omega \cap \Sigma)$. For a number $t$ suffiently small, we let $E_{th}$ be the set whose boundary $\pp E_{th} $ is parametrized as
$$
\pp E_{th} =  \{ x+ th(x) \nu(x) \ /\ x\in \pp E \},
$$
with exterior normal vector close to $\nu$.
 The first  variation of the perimeter along these normal perturbations yields
$$
\frac d{dt} Per_s( E_{th}, \Omega) \Big |_{t=0}   =  - \int_{\Sigma} H^s_{\Sigma} h ,
 $$
 where
\begin{align}
\label{1}
 H_{\Sigma}^s (p):=
{\cb \text{p.v.}}
\int_{\R^N} \frac{\chi_E (x)-\chi_{\R^N \setminus E}(x)} {|x-p|^{N+s}} \, d x
\quad \text{for } p \in \Sigma
.
\end{align}
This integral is  well-defined in the principal value sense provided that $\Sigma$ is regular near $p$.
We say that the set $\Sigma = \pp E$ is a {\em nonlocal minimal surface in $\Omega$} if the surface $\Sigma\cap \Omega$ is sufficiently regular, and it
satisfies the nonlocal minimal surface equation
\begin{equation*}
H_{\Sigma}^s (p)=0 \quad
\hbox{for all $p\in \Sigma\cap \Omega$.}
\end{equation*}
We may naturally call $H_{\Sigma}^s(p)$ the {\em nonlocal mean curvature} of $\Sigma$ at $p$.

 \medskip
 Let $\Sigma=\pp E$ be a nonlocal minimal surface. As we will prove in Section 4,
 the {\em second variation} of the $s$-perimeter in $\Omega$ can be computed for functions $h$ smooth and compactly supported in $\Sigma\cap\Omega$ as
\begin{equation}
\label{j1}
\frac {d^2}{dt^2} Per_{s} (E_{th},\Omega) \Big|_{t=0} \ =\ -2 \int_\Sigma  \JJ^s_\Sigma [h]\, h
\end{equation}
where $\JJ^s_\Sigma [h]$ is the {\em nonlocal Jacobi operator} given by

\begin{align}
\label{nonlocal jac}
 \JJ^s_\Sigma [h](p)
=
{\cb \text{p.v.}}
\int_{\Sigma}
\frac{h(x)-h(p)}{|p-x|^{N+s}} dx
+
h(p)
\int_\Sigma
\frac{\langle \nu (p)- \nu(x), \nu (p) \rangle }{|p-x|^{N+s}} dx,
\quad  p \in \Sigma.
\end{align}
In agreement with formula \eqref{j1}, we say that an  $s$-minimal surface $\Sigma$ is {\em stable} in $\Omega$ if
$$ -\int_\Sigma  \JJ^s_\Sigma [h]\, h\, \ge \, 0\foral h\in C_0^\infty(\Sigma\cap\Omega) .
$$

\medskip
A basic example of a stable nonlocal minimal surface is a {\em nonlocal area minimizing surface}.
We say that $\Sigma = \pp E$ is nonlocal area minimizing in $\Omega$ if
\begin{equation}
\label{m1}
Per_s( E,\Omega ) \le Per_s(F,\Omega)
\end{equation}
for all $F$ such that $(E\setminus F)\cup(F\setminus E)$ is compactly contained in $\Omega$.
In \cite{caffarelli-roquejoffre-savin}, Caffarelli, Roquejoffre and Savin proved that if $\Omega$  and $E_0 \subset \R^N \setminus \Omega$ are given, and sufficiently regular, then there exists a set $E$ with $E\cap (\R^N \setminus\Omega) =E_0$ which satisfies \eqref{m1}. They proved that  $\Sigma =\partial E$ is a surface of class  $C^{1,\alpha}$ outside a closed set of Hausdorff dimension $N-2$.

\medskip
 In this paper we will focus our attention on {\em nonlocal minimal cones}.
By a (solid) cone in $\R^N$, we mean a set of the form
$$
E = \{ tx\ / \ t>0, \ x\in \OO\}
$$
where $\OO$ is a regular open subset of the sphere $S^{N-1}$. The cone (mantus) $\Sigma  =\pp E$ is an $(N-1)$-dimensional surface which is regular, except at the origin.

Existence or non-existence of area minimizing cones for a given dimension is a crucial element in the classical regularity theory of minimal surfaces.
 Simons \cite{simons} proved that {\em no stable minimal cone exists in dimension $N\le 7$, except for hyperplanes}.  This result is a main ingredient in regularity theory: it implies that area minimizing surfaces must be smooth outside a closed set of Hausdorff dimension $N-8$.

\medskip
Savin and Valdinoci \cite{savin-valdinoci}, by proving the nonexistence of a nonlocal minimizing cone in $\R^2$, established the regularity of any nonlocal minimizing surface
outside a set of Hausdorff dimension $N-3$, thus improving the original result in \cite{caffarelli-roquejoffre-savin}.

In \cite{caffarelli-valdinoci}, Caffarelli and Valdinoci  proved that regularity of non-local minimizers holds up to a $(N-8)$-dimensional set, provided that $s$ is sufficiently close to $1$.

\medskip
The purpose of this paper is to analyze a specific class of nonlocal minimal cones.
Let $n,m\ge 1$, $n+m= N$ and $\A>0$. Let us call
\begin{equation}
\label{cA}
C_\alpha = \{ x = (y,z)\in \R^m\times \R^n \ /  \  |z| = \alpha |y| \ \}.
\end{equation}
It is a well-known fact that $C_\A$ is a minimal surface in $\R^N\setminus\{0\}$  (its mean curvature equals zero) if and
only $$  n\ge 2, \ m\ge 2, \quad \alpha = \sqrt{\frac{n-1}{m-1}}.$$
We call this minimal Lawson cone $C_m^n$ (\cite{lawson}). As for the stability-minimizing character of these cones, the result of Simons \cite{simons} tells us that they are all unstable for $m+n\le 7$. Simons also proved that the cone $C_4^4$ is stable and conjectured that it was minimizing.
Bombieri, De Giorgi and Giusti in \cite{bdg} found a family of disjoint minimal surfaces asymptotic to the cone, foliating $\R^4 \times \R^4$.
This implies $ \gamma= C_4^4$ is area minimizing.
For $N>8$ the cones $C_m^n$ are all area minimizing. For $N=8$ they are area minimizing if and only if $|m-n|\le 2$.
These facts were established by Lawson \cite{lawson}  and Simoes \cite{simoes},
{\cb see also \cite{miranda,concus-miranda,benarros-miranda,davini}.}

\medskip
For the non-local scenario we find the existence of analogs of the cones $C_m^n$.

\begin{theorem}
\label{thm1}
For any given $m\ge 1$, $n\ge 1$, $0<s<1$, there is a unique $\alpha=\alpha(s,m,n)>0$ such that $C_\alpha =  \{ x = (y,z)\in \R^m\times \R^n \ /  \  |z| = \alpha |y| \ \}$  is a nonlocal minimal cone. We call this cone $C_m^n (s)$.
\end{theorem}

The above result includes the existence of a minimal cone
{\cb$C_m^1(s)$, $m\geq 1$}. Such an object does not exist in the classical setting for $C_m^n$ is defined only if $n,m\ge 2$.

We have found a (computable) criterion to decide  whether or not  $C_m^n(s)$ is stable. As a consequence
we find the following result for $s$ close to $0$ which shows a sharp contrast with the classical case.

\begin{theorem}
\label{thm stability}
There is a $s_0>0$ such that  for each $s\in (0,s_0)$, all minimal cones $C_m^n(s)$ are  unstable if $N= m+ n\leq 6$ and
stable if  $N= 7$.
\end{theorem}

We recall that in the classical case $ C_m^n$ is unstable for $N=7$.
It is natural to conjecture
that the above cones for $N=7$ are minimizers of perimeter.
Being that the case, the best regularity possible for small $s$ would be
up to an $(N-7)$-dimensional set.

\medskip

As far as we know, at this moment, there are no examples of regular nontrivial nonlocal minimal surfaces (\cite{val}).  Formula (\ref{1new}) suggests that for $s$ close to $1$ there may be nontrivial nonlocal minimal surfaces close to the classical ones. In a forthcoming paper \cite{ddw} we prove that this is indeed the case. We construct  nonlocal catenoids as well as nonlocal Costa surfaces for $s$ close to $1$ by interpolating the classical minimal surfaces in compact regions with the nonlocal Lawson's cones $C_m^1$ far away. Thus these nonlocal catenoids can be considered as foliations of the nonlocal Lawson's cones $C_m^1$. A natural question, as in the classical minimal cones case (\cite{davini}), is the existence of foliations for general nonlocal Lawson's cones $C_m^n$.

\medskip

In section~\ref{sect exist unique} we prove theorem~\ref{thm1} and in section~\ref{sect s=0} we show that also for $s=0$ there is a unique minimal cone. In section~\ref{sect jacobi} we obtain formula \eqref{nonlocal jac} for the  nonlocal Jacobi operator and section~\ref{sect stability} is devoted to the proof of theorem~\ref{thm stability}.

\section{Existence and uniqueness}
\label{sect exist unique}

Let us write
\begin{align}
\label{def solid cone}
E_\alpha  = \{ x = (y,z) \, : \ y \in \R^m, \, z \in \R^n, \, |z| > \alpha |y| \ \} ,
\end{align}
so that $C_\alpha = \partial E_\alpha$
is the cone defined in \eqref{cA}.

\medskip
\noindent
{\bf Proof of theorem~\ref{thm1}.}

\noindent{\bf Existence.}
We fix $N$, $m$, $n$ with $N=m+n$, $n\leq m$ and also fix $0<s<1$.
If $m=n$ then $C_1$ is a minimal cone, since \eqref{1} is satisfied by symmetry. So we concentrate next on the case $n<m$.

Before proceeding we remark that for a cone $C_\alpha$ the quantity appearing in \eqref{1} has a fixed sign for all $p\in C_\alpha$, $p\not=0$, since by rotation we can always assume that $p = r p_\alpha$ for some $r>0$ where
\begin{align*}
p_{\alpha} = \frac{1}{\sqrt{1+\alpha^2}}(e_1^{(m)}, \alpha e_1^{(n)} )
\end{align*}
with
\begin{align}
\label{notation e}
e_1^{(m)} = (1,0,\ldots,0) \in \R^m
\end{align}
and similarly for $ e_1^{(n)}$. Then we observe that
$$
\text{p.v.}
\int_{\R^N} \frac{\chi_{E_\alpha} (x)-\chi_{E_\alpha^c}(x)} {|x-r p_{\alpha}|^{N+s}} \, d x =
\frac{1}{r^s}
\text{p.v.}
\int_{\R^N} \frac{\chi_{E_\alpha}  (x)-\chi_{E_\alpha^c} (x)} {|x-p_{\alpha}|^{N+s}} \, d x .
$$
Let us define
\begin{align}
\label{def Halpha}
H(\alpha )
=
\text{p.v.}
\int_{\R^N} \frac{\chi_{E_\alpha} (x)-\chi_{E_\alpha^c}(x)} {|x-p_{\alpha}|^{N+s}} \, d x
\end{align}
and note that it is a continuous function of $\alpha \in (0,\infty)$.

\medskip
\noindent{\bf Claim 1.} We have
\begin{align}
\label{h1}
H(1) \leq 0.
\end{align}

Indeed, write $y \in \R^m$ as $y = (y_1,y_2)$ with $y_1 \in \R^n$ and $y_2\in\R^{m-n}$.
Abbreviating $e_1 = e_1^{(n)} = (1,0,\ldots,0) \in \R^n$ we rewrite
\begin{align*}
H(1)
&= \lim_{\delta\to0}
\int_{\R^N \setminus B(p_{1},\delta)} \frac{\chi_{E_1} (x)-\chi_{E_1^c}(x)} {|x-p_{1}|^{N+s}} \, d x
\\
&
=
\lim_{\delta\to0}
\int_{A_\delta} \frac{1}{( |y_1 - \frac1{\sqrt2}e_1|^2 + |y_2|^2 + |z-\frac1{\sqrt2}e_1|^2)^{\frac{N+s}{2}} }
\\
&\qquad -
 \lim_{\delta\to0}
\int_{B_\delta} \frac{1}{( |y_1 - \frac1{\sqrt2}e_1|^2 + |y_2|^2 + |z-\frac1{\sqrt2}e_1|^2)^{\frac{N+s}{2}} } ,
\end{align*}
where
\begin{align*}
A_\delta
& = \{ |z|^2>|y_1 |^2 + |y_2|^2, \
|y_1 - \frac1{\sqrt2}e_1|^2 + |y_2|^2 + |z-\frac1{\sqrt2}e_1|^2>\delta^2\}
\\
B_\delta
&= \{ |z|^2<|y_1 |^2 + |y_2|^2,  \
|y_1 - \frac1{\sqrt2}e_1|^2 + |y_2|^2 + |z-\frac1{\sqrt2}e_1|^2>\delta^2 \} .
\end{align*}
But the first integral can be rewritten as
\begin{align*}
& \int_{A_\delta} \frac{1}{( |y_1 - \frac1{\sqrt2}e_1|^2 + |y_2|^2 + |z-\frac1{\sqrt2}e_1|^2)^{\frac{N+s}{2}} }
\\
& =
\int_{\tilde A_\delta} \frac{1}{( |y_1 - \frac1{\sqrt2}e_1|^2 + |y_2|^2 + |z-\frac1{\sqrt2}e_1|^2)^{\frac{N+s}{2}} }
\end{align*}
where
$$
\tilde A_\delta = \{
|y_1|^2>|z |^2 + |y_2|^2, \
|y_1 - \frac1{\sqrt2}e_1|^2 + |y_2|^2 + |z-\frac1{\sqrt2}e_1|^2>\delta^2\}
$$
(we just have exchanged $y_1$ by $z$ and noted that the integrand is symmetric in these variables).
But $\tilde A_\delta \subset B_\delta$ and so
\begin{align*}
& \int_{\R^N \setminus B(p_{1},\delta)} \frac{\chi_{E_1} (x)-\chi_{E_1^c}(x)} {|x-p_{1}|^{N+s}} \, d x
\\
& =
- \int_{B_\delta \setminus \tilde A_\delta}
\frac{1}{( |y_1 - \frac1{\sqrt2}e_1|^2 + |y_2|^2 + |z-\frac1{\sqrt2}e_1|^2)^{\frac{N+s}{2}} }
\leq 0.
\end{align*}
This shows the validity of \eqref{h1}.

\medskip
\noindent{\bf Claim 2.} We have
\begin{align}
\label{limit H 0}
H(\alpha) \to+\infty \quad\text{as }\alpha\to 0.
\end{align}

Let $0<\delta<1/2$ be fixed and write
$$
H(\alpha) = I_\alpha + J_\alpha
$$
where
\begin{align*}
I_\alpha
&=\int_{\R^N \setminus B(p_{\alpha},\delta)} \frac{\chi_{E_\alpha} (x)-\chi_{E_\alpha^c}(x)} {|x-p_{\alpha}|^{N+s}} \, d x
\\
J_\alpha &=
\text{p.v.} \int_{B(p_{\alpha},\delta)} \frac{\chi_{E_\alpha} (x)-\chi_{E_\alpha^c}(x)} {|x-p_{\alpha}|^{N+s}} \, d x .
\end{align*}
With $\delta$ fixed
\begin{align}
\label{limit I}
\lim_{\alpha\to 0} I_\alpha
= \int_{\R^N \setminus B(p_{\alpha},\delta)} \frac{1} {|x-p_{0}|^{N+s}} \, d x  >0.
\end{align}
For $J_\alpha$ we make a change of variables $x = \alpha \tilde x + p_\alpha$ and obtain
\begin{align}
\label{scaling J}
J_\alpha
=
\text{p.v.} \int_{B(p_{\alpha},\delta)} \frac{\chi_{E_\alpha} (x)-\chi_{E_\alpha^c}(x)} {|x-p_{\alpha}|^{N+s}} \, d x
=
\frac{1}{\alpha^s}
\text{p.v.}
\int_{B(0,\delta/\alpha)}
\frac{\chi_{F_\alpha}(\tilde x) -\chi_{F_\alpha^c} (\tilde x)}{|\tilde x|^{N+s}} d \tilde x
\end{align}
where $F_\alpha = \frac1\alpha ( E_\alpha-p_\alpha)$. But
$$
\text{p.v.}
\int_{B(0,\delta/\alpha)}
\frac{\chi_{F_\alpha}(\tilde x) -\chi_{F_\alpha^c} (\tilde x)}{|\tilde x|^{N+s}} d \tilde x
\to
\text{p.v}
\int_{\R^N}
\frac{\chi_{F_0}(x) -\chi_{F_0^c} (x)}{|x|^{N+s}} dx
$$
as $\alpha\to 0$ where
$F_0 = \{ x=(y,z) : y\in \R^m, z\in \R^n, \ |z+e_1^{(n)}|>1 \} $.
But writing $z=(z_1,\ldots,z_n)$ we see that
\begin{align*}
\text{p.v}
\int_{\R^N}
\frac{\chi_{F_0}(x) -\chi_{F_0^c} (x)}{|x|^{N+s}} d x
&\geq
\text{p.v}
\int_{\R^N}
\frac{\chi_{[z_1>0 \text{ or } z_1<-2]} -\chi_{[-2<z_1<0]}}{|x|^{N+s}} dx
\\
&\geq
\int_{\R^N}
\frac{\chi_{ [ \ |z_1| >2 \ ]} }{|x|^{N+s}}  d x
\end{align*}
and this number is positive. This and \eqref{scaling J} show that $J_\alpha\to+\infty$ as $\alpha\to0$ and combined with \eqref{limit I} we obtain the desired conclusion.


By  \eqref{h1}, \eqref{limit H 0} and continuity we obtain the existence of $\alpha \in (0,1]$ such that $H(\alpha)=0$.

\bigskip
\noindent
{\bf Uniqueness.}
Consider 2 cones $C_{\alpha_1}$, $C_{\alpha_2}$ with $\alpha_1>\alpha_2>0$, associated to solid cones $E_{\alpha_1}$ and $E_{\alpha_2}$.
We claim that there is a rotation $R$ so that $R(E_{\alpha_1}) \subset E_{\alpha_2}$ (strictly) and that
$$
H(\alpha_1) =
\text{p.v.}
\int_{\R^N}
\int_{\R^N} \frac{\chi_{R(E_{\alpha_1})} (x)-\chi_{R(E_{\alpha_1})^c}(x)} {|x-p_{\alpha_2}|^{N+s}} \, d x .
$$
Note that the denominator in the integrand is the same that appears in \eqref{def Halpha} for $\alpha_2$ and then
\begin{align}
\nonumber
H(\alpha_1)
&=
\text{p.v.}
\int_{\R^N}
\int_{\R^N} \frac{\chi_{R(E_{\alpha_1})} (x)-\chi_{R(E_{\alpha_1})^c}(x)} {|x-p_{\alpha_2}|^{N+s}} \, d x
\\
\label{strict ineq}
&<
\text{p.v.}
\int_{\R^N}
\int_{\R^N} \frac{\chi_{E_{\alpha_2}} (x)-\chi_{E_{\alpha_2}^c}(x)} {|x-p_{\alpha_2}|^{N+s}} \, d x = H(\alpha_2) .
\end{align}
This shows that $H(\alpha)$ is decreasing in $\alpha$ and hence the uniqueness.
To construct the rotation let us write as before $x = (y,z) \in \R^N$, with $y\in\R^m$, $z\in \R^n$, and $y = (y_1,y_2)$ with $y_1\in\R^n$, $y_2 \in \R^{m-n}$ (we assume alway $n\leq m$). Let us write the vector $(y_1,z)$ in spherical coordinates of $\R^{2n}$ as follows
$$
y_1
= \rho
\left[
\begin{matrix}
\cos(\varphi_1)\\
\sin(\varphi_1) \cos(\varphi_2)\\
\sin(\varphi_1) \sin(\varphi_2) \cos(\varphi_3)\\
\vdots\\
\sin(\varphi_1) \sin(\varphi_2) \sin(\varphi_3)\ldots \sin(\varphi_{n-1}) \cos(\varphi_n)\\
\end{matrix}
\right]
$$
$$
z
= \rho
\left[
\begin{matrix}
\sin(\varphi_1) \sin(\varphi_2) \sin(\varphi_3)\ldots \sin(\varphi_{n}) \cos(\varphi_{n+1})\\
\vdots\\
\sin(\varphi_1) \sin(\varphi_2) \sin(\varphi_3)\ldots \sin(\varphi_{2n-2}) \cos(\varphi_{2n-1})\\
\sin(\varphi_1) \sin(\varphi_2) \sin(\varphi_3)\ldots \sin(\varphi_{2n-2}) \sin(\varphi_{2n-1})\\
\end{matrix}
\right]
$$
where $\rho>0$, $\varphi_{2n-1}\in[0,2\pi)$, $\varphi_j\in[0,\pi]$ for $j=1,\ldots,2n-2$. Then
$$
|z|^2 = \rho^2 \sin(\varphi_1)^2 \sin(\varphi_2)^2 \ldots \sin(\varphi_n)^2
,
\quad |y_1|^2+|z|^2 = \rho^2.
$$
The equation for the solid cone $E_{\alpha_i}$, namely $|z|>\alpha_i |y|$, can be rewritten as
$$
\rho^2 \sin(\varphi_1)^2 \sin(\varphi_2)^2 \ldots \sin(\varphi_n)^2
> \alpha_i^2 (|y_1|^2 + |y_2|^2).
$$
Adding $\alpha_i^2 |z|^2$ to both sides this is equivalent to
$$
\sin(\varphi_1)^2 \sin(\varphi_2)^2 \ldots \sin(\varphi_n)^2
> \sin(\beta_i)^2
( 1+\frac{ |y_2|^2}{\rho^2})
$$
where $\beta_i = \arctan(\alpha_i)$.
We let $\theta = \beta_1 - \beta_2 \in (0,\pi/2)$,
and define the rotated cone $R_\theta(E_{\alpha_1})$ by the equation
$$
\sin(\varphi_1+\theta)^2 \sin(\varphi_2)^2 \ldots \sin(\varphi_n)^2
> \sin(\beta_1)^2
( 1+\frac{ |y_2|^2}{\rho^2}) .
$$
We want to show that $R_\theta(E_{\alpha_1}) \subset E_{\alpha_2}$. To do so, it suffices to prove that for any given $t\geq 1$, if $\varphi$ satisfies the inequality $ |\sin(\varphi+\theta)|> \sin(\beta_1) t $
then it also satisfies $ |\sin(\varphi)|> \sin(\beta_2) t $.  This in turn can be proved from the inequality
\begin{align*}
\arccos(\sin(\beta_1) t) + \theta < \arccos(\sin(\beta_2) t)
\end{align*}
for $ 1< t \leq \frac{1}{\sin(\beta_1)}$. For $t=1$ we have equality by definition of $\theta$. The inequality for $ 1< t \leq \frac{1}{\sin(\beta_1)}$  can be checked by computing a derivative with respect to $t$. The strict inequality in \eqref{strict ineq} is because $R(E_{\alpha_1}) \subset E_{\alpha_2}$ strictly.
\qed


\section{Minimal cones for $s=0$}
\label{sect s=0}

In this section we derive the limiting value $\alpha_0 = \lim_{s\to0} \alpha_s$ where $\alpha_s$ is such that $C_{\alpha_s}$ is an $s$-minimal cone.

\begin{prop}
\label{prop alpha s=0}
Assume that $n\leq m$ in \eqref{def solid cone}, $N=m+n$.
The number $\alpha_0$ is the unique solution to
$$
\int_{\alpha}^\infty
\frac{t^{n-1}}
{(1 + t^2  )^\frac{N}2}
d t - \int_0^{\alpha}
\frac{t^{n-1}}
{(1 + t^2  )^\frac{N}2}
d t
=0 .
$$
\end{prop}
\noindent{\bf Proof.}
We write $x = (y,z)\in\R^N$ with $y\in \R^m$, $z\in \R^n$.
Let us assume in the rest of the proof that $n\geq 2$. The case $n=1$ is similar.
We evaluate the integral in \eqref{1} for the point $p=(e_1^{(m)},\alpha e_1^{(n)})$
using spherical coordinates for $y = r \omega_1$ and $z = \rho \omega_2$
where $r,\rho>0$ and
\begin{align}
\label{omega1}
\omega_1 = \left[
\begin{matrix}
\cos(\theta_1)\\
\sin(\theta_1) \cos(\theta_2)\\
\vdots\\
\sin(\theta_1) \sin(\theta_2) \ldots \sin(\theta_{m-2}) \cos(\theta_{m-1})\\
\sin(\theta_1) \sin(\theta_2) \ldots
\sin(\theta_{m-2}) \sin(\theta_{m-1})
\end{matrix}
\right]
\end{align}
\begin{align}
\label{omega2}
\omega_2 = \left[
\begin{matrix}
\cos(\varphi_1)\\
\sin(\varphi_1) \cos(\varphi_2)\\
\vdots\\
\sin(\varphi_1) \sin(\varphi_2) \ldots \sin(\varphi_{n-2}) \cos(\varphi_{n-1})\\
\sin(\varphi_1) \sin(\varphi_2) \ldots
\sin(\varphi_{n-2}) \sin(\varphi_{n-1})
\end{matrix}
\right] ,
\end{align}
where $\theta_j\in [0,\pi]$ for $j=1,\ldots,m-2$, $\theta_{m-1} \in [0,2\pi]$, $\varphi_j\in [0,\pi]$ for $j=1,\ldots,n-2$, $\varphi_{n-1} \in [0,2\pi]$.
Then
$$
|(y,z)-(e_1^{(m)},\alpha e_1^{(n)})|^2
=r^2 + 1 - 2 r \cos(\theta_1) + \rho^2 + \alpha^2 -2\rho\alpha\cos(\varphi_1) .
$$
Assuming that $\alpha = \alpha_s>0$ is such that $C_{\alpha_s}$ is an $s$-minimal cone, \eqref{1} yields the following equation for $\alpha$
\begin{align}
\label{eq alpha}
\text{p.v.}
\int_0^\infty
r^{m-1}
( A_{\alpha,s}(r) - B_{\alpha,s}(r) ) d r = 0
\end{align}
where
\begin{align*}
A_{\alpha,s}(r)
&=
\int_{r\alpha}^\infty
\int_0^\pi
\int_0^\pi
\frac{\rho^{n-1}
\sin(\theta_1)^{m-2}
\sin(\varphi_1)^{n-2}
}{(r^2 + 1 - 2 r \cos(\theta_1) + \rho^2 + \alpha^2 -2\rho\alpha\cos(\varphi_1) )^\frac{N+s}2}
d\theta_1
d\varphi_1
d\rho
\\
 B_{\alpha,s}(r)
 &=
\int_0^{r\alpha}
\int_0^\pi
\int_0^\pi
\frac{\rho^{n-1}
\sin(\theta_1)^{m-2}
\sin(\varphi_1)^{n-2}
}{(r^2 + 1 - 2 r \cos(\theta_1) + \rho^2 + \alpha^2 -2\rho\alpha\cos(\varphi_1) )^\frac{N+s}2}
d\theta_1
d\varphi_1
d\rho ,
\end{align*}
which are well defined for $r\not=1$.
Setting $\rho = r t$ we get
\begin{align*}
& A_{\alpha,s}(r)
\\
& =
r^{-m-s}
\int_{\alpha}^\infty
\int_0^\pi
\int_0^\pi
\frac{t^{n-1}\sin(\varphi_1)^{m-2}
\sin(\theta_1)^{n-2}}
{(1 + \frac1{r^2} - \frac2r  \cos(\theta_1) + t^2 + \frac{\alpha^2}{r^2} - \frac2r t \alpha\cos(\varphi_1) )^\frac{N+s}2}
d\theta_1
d\varphi_1
d t
\\
& =
c_{m,n}
r^{-m-s}
\int_{\alpha}^\infty
\frac{t^{n-1}}
{(1 + t^2  )^\frac{N+s}2}
d t + O(r^{-m-s-1})
\end{align*}
as $r\to\infty$ and this is uniform in $s$ for $s>0$ small.
Here $c_{m,n}>0$ is some constant.
Similarly
\begin{align*}
B_{\alpha,s}(r) &=
c_{m,n}
r^{-m-s}
\int_0^{\alpha}
\frac{t^{n-1}}
{(1 + t^2  )^\frac{N+s}2}
d t + O(r^{-m-s-1})
\end{align*}
Then \eqref{eq alpha} takes the form
\begin{align*}
0&=\int_0^2 \ldots dr + \int_2^\infty \ldots dr
= O(1) + C_s(\alpha) \int_2^\infty r^{-1-s} dr
=
O(1) + \frac{2^{-s}}s C_s(\alpha)
\end{align*}
where
$$
C_s(\alpha) = \int_{\alpha}^\infty
\frac{t^{n-1}}
{(1 + t^2  )^\frac{N+s}2}
d t - \int_0^{\alpha}
\frac{t^{n-1}}
{(1 + t^2  )^\frac{N+s}2}
d t
$$
and $O(1)$ is uniform as $s\to0$, because $0<\alpha_s \leq 1$ by theorem~\ref{thm1},
and the only singularity in \eqref{eq alpha} occurs at $r = 1$.
This implies that $\alpha_0 = \lim_{s\to0} \alpha_s$ has to satisfy $C_0(\alpha_0)=0$.
\qed

\section{The Jacobi operator}
\label{sect jacobi}

In this section we prove formula \eqref{j1} and derive the formula for the nonlocal Jacobi operator \eqref{nonlocal jac}.

Let $E \subset \R^N$ be an open set with smooth boundary and $\Omega$ be a bounded open set.
Let $\nu$ be the unit normal vector field  of $\Sigma = \pp E$ pointing to the exterior of $E$.
Given $h\in C_0^\infty (\Omega \cap \Sigma)$ and $t$  small, let $E_{th}$ be the set whose boundary $\pp E_{th} $ is parametrized as
$$
\pp E_{th} =  \{ x+ th(x) \nu(x) \ /\ x\in \pp E \},
$$
with exterior normal vector close to $\nu$.

\begin{prop}
\label{prop second var}
For $h\in C_0^\infty (\Omega \cap \Sigma)$
\begin{align}
\label{second var}
\frac {d^2}{dt^2} Per_{s} (E_{th},\Omega) \Big|_{t=0} \ =\ -2 \int_\Sigma  \JJ^s_\Sigma [h]\, h
-\int_{\Sigma} h^2 H H_\Sigma^s ,
\end{align}
where $ \JJ^s_\Sigma $ is the nonlocal Jacobi operator defined in \eqref{nonlocal jac}, $H$ is the classical mean curvature of $\Sigma$ and $H_\Sigma^s$ is the nonlocal mean curvature defined in \eqref{1}.
\end{prop}

In case that $\Sigma$ is a nonlocal minimal surface in $\Omega$ we obtain formula \eqref{j1}.
Another related formula is the following.

\begin{prop}
\label{prop derivative H}
Let $\Sigma_{th} = \partial E_{th}$.
For $p\in \Sigma$ fixed let
$p_t = p + t h(p) \nu(p) \in \Sigma_{th}$.
Then for $h\in C^\infty ( \Sigma) \cap L^\infty(\Sigma)$
\begin{align}
\label{der mean C}
\frac{d}{d t}
H_{\Sigma_{th}}^s(p_t)
\Big|_{t=0}
= 2 \JJ_\Sigma^s[h] (p) .
\end{align}
\end{prop}

A consequence of proposition~\ref{prop derivative H} is that entire nonlocal minimal graphs are stable.

\begin{corollary}
\label{coro entire stable}
Suppose that $\Sigma = \partial E$ with
$$
E = \{ (x',F(x')) \in \R^N: x'\in \R^{N-1} \}
$$
is a nonlocal minimal surface. Then
\begin{align}
\label{linearly stable}
-\int_\Sigma  \JJ^s_\Sigma [h]\, h\, \ge \, 0\foral h\in C_0^\infty(\Sigma).
\end{align}
\end{corollary}

\medskip
\noindent
{\bf Proof of proposition~\ref{prop second var}.}
Let
$$
K_\delta (z) = \frac{1}{|z|^{N+s}} \eta_\delta(z)
$$
where
$\eta_\delta(x) = \eta(x/\delta)$
($\delta >0$) and $\eta \in C^\infty(\R^N)$ is a radially symmetric cut-off function with $\eta(x)=1 $ for $|x|\geq 2$, $\eta(x) = 0$ for $|x|\leq 1$.

Consider
\begin{align}
\label{def per delta}
Per_{s,\delta}(E_{th},\Omega)
=
\int_{E_{th}\cap\Omega} \int_{\R^N \setminus E_{th}}
K_\delta(x-y) \, dy \, d x
+  \int_{E_{th}\setminus\Omega} \int_{\Omega \setminus E_{th}}
K_\delta(x-y) d y d x .
\end{align}
We will show that $\frac{d^2}{dt^2}
Per_{s,\delta}(E_{th},\Omega)  $ approaches a certain limit $D_2(t)$ as $\delta \to 0$, uniformly for $t$ in a neighborhood of $0$ and that
\begin{align*}
D_2(0)&=
-2 \int_\Sigma  \JJ^s_\Sigma [h]\, h
-\int_{\Sigma} h^2 H H_\Sigma^s .
\end{align*}

First we need some extensions of $\nu$ and $h$ to $\R^N$.
To define them, let $K \subset \Sigma$ be the support of $h$ and
$U_0$ be an open bounded neighborhood of $K$ such that  for any $x\in U_0$, the closest point $\hat x \in \Sigma$ to $x$ is unique and defines a smooth function of $x$.
We also take $U_0$ smaller if necessary as to have $\overline U_0 \subset \Omega$.
Let $\tilde \nu:\R^N \to \R^N$ be a globally defined smooth unit vector field such that
$\tilde \nu(x)=\nu(\hat x)$ for $x\in U_0$.
We also extend $h$ to $\tilde h: \R^N \to \R$ such that it is smooth with compact support contained in $\Omega$ and $\tilde h(x) = h(\hat x)$ for $x\in U_0$.
From now one we omit the tildes ($\tilde\ $) in the definitions of the extensions of $\nu$ and $h$.
For $t$ small  $\bar x \mapsto \bar x + t h(\bar x) \nu(\bar x) $ is a global diffeomorphism in $\R^N$. Let us write
$$
u(\bar x) =  h(\bar x) \nu(\bar x)  \quad \text{for } \bar x\in \R^N ,
$$
$$
\nu = (\nu^1,\ldots,\nu^N), \quad u=  (u^1,\ldots,u^N)
$$
and let
$$
J_t(\bar x) = J_{id + t u}(\bar x)
$$
be the Jacobian determinant of  $id + t u$.

We change variables
$$
x = \bar x + t u (\bar x) ,
\quad
y = \bar y + t u (\bar y) ,
$$
in \eqref{def per delta}
\begin{align*}
Per_{s,\delta}(E_{th},\Omega)
& =
\int_{E\cap \phi_t(\Omega) } \int_{\R^N \setminus E}
K_\delta(x-y)
J_t(\bar x) J_t(\bar y)
d\bar yd \bar x,
\\
& \qquad
+
\int_{E\setminus \phi_t(\Omega)} \int_{\phi_t(\Omega) \setminus E}
K_\delta(x-y)
J_t(\bar y) d\bar yd \bar x ,
\end{align*}
where $\phi_t$ is the inverse of the map $\bar x \mapsto \bar x + t u(\bar x)$.

Differentiating with respect to $t$:
\begin{align*}
\frac{d}{dt}
Per_{s,\delta}(E_{th},\Omega)
 & =
 \int_{E\cap \phi_t(\Omega) } \int_{\R^N \setminus E}
 \Big[
\nabla K_\delta(x-y) (u(\bar x) - u(\bar y)) J_t(\bar x) J_t(\bar y)
\\
& \qquad
+
K_\delta(x-y) ( J_t'(\bar x) J_t(\bar y) + J_t(\bar x) J_t'(\bar y) )
\Big]
d\bar y d \bar x
\\
& \qquad
+
\int_{E\setminus \phi_t(\Omega)} \int_{\phi_t(\Omega) \setminus E}
\Big[
\nabla K_\delta(x-y) (u(\bar x) - u(\bar y)) J_t(\bar x) J_t(\bar y)
\\
&\qquad
+
K_\delta(x-y) ( J_t'(\bar x) J_t(\bar y) + J_t(\bar x) J_t'(\bar y) )
\Big]
d\bar y d \bar x ,
\end{align*}
where
$$
J_t'(\bar x ) = \frac{d}{dt} J_t(\bar x) .
$$
Note that there are no integrals on $\partial \phi_t(\Omega)$ for $t$ small because $u$ vanishes in a neighborhood of $\partial \Omega$.

Since the integrands in $\frac{d}{dt}
Per_{s,\delta}(E_{th},\Omega) $ have compact support contained in $\phi_t(\Omega)$ ($t$ small), we can write
\begin{align*}
\frac{d}{dt}
Per_{s,\delta}(E_{th},\Omega)
&=
\int_{E} \int_{\R^N\setminus E}
\Big[
\nabla K_\delta(x-y) (u(\bar x) - u(\bar y)) J_t(\bar x) J_t(\bar y)
\\
& \qquad +
 K_\delta(x-y) ( J_t'(\bar x) J_t(\bar y) + J_t(\bar x) J_t'(\bar y) )
\Big]
d\bar y d \bar x .
\end{align*}
Differentiating once more
$$
\frac{d^2}{dt^2}
Per_{s,\delta}(E_{th},\Omega)
= A(\delta,t) + B(\delta,t) + C(\delta,t)
$$
where
\begin{align*}
A(\delta,t) & =
\int_{E} \int_{\R^N\setminus E}
D^2 K_\delta(x-y) (u(\bar x) - u(\bar y)) (u(\bar x) - u(\bar y)) J_t(\bar x) J_t(\bar y) d\bar y d \bar x
\\
B(\delta,t) & =
2 \int_{E} \int_{\R^N\setminus E}
\nabla K_\delta(x-y) (u(\bar x) - u(\bar y))
( J_t'(\bar x) J_t(\bar y) + J_t(\bar x) J_t'(\bar y) )
d\bar y d \bar x
\\
C(\delta,t) & =
\int_{E} \int_{\R^N\setminus E}
K_\delta(x-y)
(
J_t''(\bar x) J_t(\bar y)
+2 J_t'(\bar x) J_t'(\bar y)
+ J_t(\bar x) J_t''(\bar y) )
d\bar y d \bar x .
\end{align*}

We claim that $A(\delta,t)$, $B(\delta,t)$ and $C(\delta,t)$ converge as $\delta \to 0$ for uniformly for $t$ near 0, to limit expressions $A(0,t)$, $B(0,t)$ and $C(0,t)$, which are the same as above replacing $\delta$ by 0, and that the integrals appearing in $A(0,t)$, $B(0,t)$ and $C(0,t)$ are well defined.
Indeed, we can estimate
$$
|A(\delta,t)-A(0,t)|
\leq C
\int_{x\in E \cap K_0}
\int_{ y \in E^c, |x-y|\leq 2 \delta}
\frac{1}{|x-y|^{N+s}} \, d y \, d x,
$$
where $K_0$ is a fixed bounded set.
For $x\in E \cap K_0$ we see that
$$
\int_{ y \in E^c, |x-y|\leq 2 \delta}
\frac{1}{|x-y|^{N+s}} \, d y
\leq \frac{C}{dist(x,E^c)^s},
$$
and therefore
$$
|A(\delta,t)-A(0,t)|
\leq C
\leq C
\int_{x\in E \cap K_0, \ dist(x,E^c) \leq 2 \delta}
\frac{1}{dist(x,E^c)^s} \, d x
\leq C \delta^{1-s}.
$$
The differences  $B(\delta,t)-B(0,t)$, $C(\delta,t)-C(0,t)$ can be estimated similarly.
This shows that
$$
\frac{d^2}{dt^2}
Per_{s}(E_{th},\Omega)
\Big|_{t=0}
=
\lim_{\delta\to 0}
\frac{d^2}{dt^2}
Per_{s,\delta}(E_{th},\Omega)
\Big|_{t=0}
= \lim_{\delta \to 0}
A(\delta,0)+B(\delta,0)+C(\delta,0) .
$$
In what follows we will evaluate $A(\delta,0)+B(\delta,0)+C(\delta,0)$.
At $t=0$  we have
\begin{align*}
A(\delta,0)
&=
\int_{E}\int_{\R^N\setminus E}
D_{x_i x_j}
K_\delta(x-y)(u^i(x)-u^i(y))(u^j(x)-u^j(y))
\, d y \, d x
\\
&=A_{11} + A_{12} + A_{21} + A_{22}
\end{align*}
where
\begin{align*}
A_{11} &=
\int_{E} \int_{\R^N\setminus E}
D_{x_i x_j} K_\delta(x-y) u^i(x) u^j(x)
\, d y \, d x
\\
A_{12}
&=
-\int_{E} \int_{\R^N\setminus E}
D_{x_i x_j} K_\delta(x-y) u^i(x) u^j(y)
\, d y \, d x
\\
A_{21}
&=
-\int_{E} \int_{\R^N\setminus E}
D_{x_i x_j} K_\delta(x-y) u^i(y) u^j(x)
\, d y \, d x
\\
A_{22}
&=
\int_{E} \int_{\R^N\setminus E}
D_{x_i x_j} K_\delta(x-y) u^i(y) u^j(y)
\, d y \, d x .
\end{align*}
Let us also write
\begin{align*}
B(\delta,0)
&=
2 \int_{E} \int_{\R^N\setminus E}
D_{x_j} K_\delta(x-y) (u^j(x) - u^j(y)) ( {\rm div}(u)(x) +  {\rm div}(u)(y))
\, d y \, d x
\\
&= B_{11} + B_{12} + B_{21} + B_{22} ,
\end{align*}
where
\begin{align*}
B_{11} &=
2 \int_{E} \int_{\R^N\setminus E}
D_{x_j} K_\delta(x-y) u^j(x)  {\rm div}(u)(x)
\, d y \, d x
\\
B_{12}
&=
2 \int_{E} \int_{\R^N\setminus E}
D_{x_j} K_\delta(x-y)u^j(x) {\rm div}(u)(y)
\, d y \, d x
\\
B_{21}
&=
-
2 \int_{E} \int_{\R^N\setminus E}
D_{x_j} K_\delta(x-y)u^j(y) {\rm div}(u)(x)
\, d y \, d x\\
B_{22} &=
2 \int_{E} \int_{\R^N\setminus E}
D_{y_j} K_\delta(x-y) u^j(y)  {\rm div}(u)(y)
\, d y \, d x ,
\end{align*}
and
\begin{align*}
C(\delta,0)=C_1 + C_2 + C_3 ,
\end{align*}
where
\begin{align*}
C_1 & =
\int_{E} \int_{\R^N\setminus E}
K_\delta(x-y)
\Big[ div(u)(x)^2 - tr(D u(x)^2) \Big]
\, d y \, d x
\\
C_2
&=
\int_{E} \int_{\R^N\setminus E}
K_\delta(x-y)
\Big[
div(u)(y)^2 - tr(D u(y)^2)
\Big]
\, d y \, d x
\\
C_3 &=
2 \int_{E} \int_{\R^N\setminus E}
K_\delta(x-y)
div(u)(x) div(u)(y)
\, d y \, d x .
\end{align*}

We compute
\begin{align*}
A_{11}
&=
\int_{E} \int_{\R^N\setminus E}
D_{x_i} \Big[ D_{x_j} K_\delta(x-y) u^i(x) u^j(x) \Big]
\, d y \, d x
\\
& \qquad
-
\int_{E} \int_{\R^N\setminus E}
D_{x_j} K_\delta(x-y) D_{x_i} \Big[  u^i(x) u^j(x) \Big]
\, d y \, d x
\\
&=
\int_{\partial E} \int_{\R^N\setminus E}
D_{x_j} K_\delta(x-y) u^i(x) u^j(x) \nu^i(x)
\, d y \, d x
\\
& \qquad
-
\int_{E} \int_{\R^N\setminus E}
D_{x_j} K_\delta(x-y)  \Big[  D_{x_i} u^i(x) u^j(x) +  u^i(x) D_{x_i}  u^j(x)\Big]
\, d y \, d x .
\end{align*}
Therefore
\begin{align*}
A_{11} + B_{11}
&=
\int_{\partial E} \int_{\R^N\setminus E}
D_{x_j} K_\delta(x-y) u^i(x) u^j(x) \nu^i(x)
\, d y \, d x
\\
&\qquad
+
\int_{E} \int_{\R^N\setminus E}
D_{x_j} K_\delta(x-y)  \Big[  D_{x_i} u^i(x) u^j(x) -  u^i(x) D_{x_i}  u^j(x)\Big]
\, d y \, d x .
\end{align*}
We express the first term as
\begin{align*}
\int_{\partial E} \int_{\R^N\setminus E}
&
D_{x_j} K_\delta(x-y) u^i(x) u^j(x) \nu^i(x)
\, d y \, d x
\\
&
=
- \int_{\partial E} \int_{\R^N\setminus E}
D_{y_j} K_\delta(x-y) u^i(x) u^j(x) \nu^i(x)
\, d y \, d x
\\
&=
\int_{\partial E} \int_{\partial E}
K_\delta(x-y) u^i(x) u^j(x) \nu^i(x) \nu^j(y)
\, d y \, d x
\\
&=
\int_{\partial E} \int_{\partial E}
K_\delta(x-y) h(x)^2 \nu(x) \nu(y)
\, d y \, d x .
\end{align*}
For the second term of $A_{11}+B_{11}$ let us write
\begin{align*}
\int_{E} \int_{\R^N\setminus E}
& D_{x_j} K_\delta(x-y)   D_{x_i} u^i(x) u^j(x)
\, d y \, d x
\\
& =
\int_{E} \int_{\R^N\setminus E}
D_{x_j} \Big[ K_\delta(x-y)   D_{x_i} u^i(x) u^j(x) \Big]
\, d y \, d x
\\
& \qquad -
\int_{E} \int_{\R^N\setminus E}
K_\delta(x-y) D_{x_j} \Big[   D_{x_i} u^i(x) u^j(x) \Big]
\, d y \, d x
\\
& =
\int_{\partial E} \int_{\R^N\setminus E}
K_\delta(x-y)   D_{x_i} u^i(x) u^j(x) \nu^j(x)
\, d y \, d x
\\
& \qquad -
\int_{E} \int_{\R^N\setminus E}
K_\delta(x-y) \Big[   D_{x_j x_i} u^i(x) u^j(x) + div(u)(x)^2 \Big]
\, d y \, d x .
\end{align*}
The third term of $A_{11}+B_{11}$ is
\begin{align*}
- \int_{E} \int_{\R^N\setminus E}
&
D_{x_j} K_\delta(x-y)    u^i(x) D_{x_i}  u^j(x)
\, d y \, d x
\\
&=
- \int_{E} \int_{\R^N\setminus E}
D_{x_j} \Big[ K_\delta(x-y)    u^i(x) D_{x_i}  u^j(x) \Big]
\, d y \, d x
\\
&\qquad
+
\int_{E} \int_{\R^N\setminus E}
K_\delta(x-y)    D_{x_j} \Big[ u^i(x) D_{x_i}  u^j(x) \Big]
\, d y \, d x
\\
&=
- \int_{\partial E} \int_{\R^N\setminus E}
K_\delta(x-y)    u^i(x) D_{x_i}  u^j(x) \nu^j(x)
\, d y \, d x
\\
&\qquad
+
\int_{E} \int_{\R^N\setminus E}
K_\delta(x-y)     \Big[ D_{x_j} u^i(x) D_{x_i}  u^j(x) + u^i(x) D_{x_j x_i}  u^j(x)\Big]
\, d y \, d x  .
\end{align*}
Therefore
\begin{align*}
A_{11}+B_{11}
&=
\int_{\partial E} \int_{\partial E}
K_\delta(x-y) h(x)^2 \nu(x) \nu(y)
\, d y \, d x
\\
& \qquad
+
\int_{\partial E} \int_{\R^N\setminus E}
K_\delta(x-y)
\Big[ D_{x_i} u^i(x) u^j(x) \nu^j(x)
- u^i(x) D_{x_i}  u^j(x) \nu^j(x) \Big]
\, d y \, d x
\\
& \qquad
+
\int_{E} \int_{\R^N\setminus E}
K_\delta(x-y) \Big[  D_{x_j} u^i(x) D_{x_i}  u^j(x)-    {\rm div}(u)(x)^2 \Big]
\, d y \, d x ,
\end{align*}
so that
\begin{align*}
A_{11} + B_{11} + C_1
&=
\int_{\partial E} \int_{\partial E}
K_\delta(x-y) h(x)^2 \nu(x) \nu(y)
\, d y \, d x
\\
& \qquad
+
\int_{\partial E} \int_{\R^N\setminus E}
K_\delta(x-y)
\Big[ D_{x_i} u^i(x) u^j(x) \nu^j(x)
- u^i(x) D_{x_i}  u^j(x) \nu^j(x) \Big]
\, d y \, d x .
\end{align*}
But using $u = \nu h$ and $div(\nu) = H$ where $H$ is the mean curvature of $\partial E$ we have
\begin{align*}
D_{x_i} u^i(x) u^j(x) \nu^j(x)
- u^i(x) D_{x_i}  u^j(x) \nu^j(x)
&=
h(x)^2 H(x)
\end{align*}
and therefore
\begin{align*}
A_{11} + B_{11} + C_1
&=
\int_{\partial E} \int_{\partial E}
K_\delta(x-y) h(x)^2 \nu(x) \nu(y)
\, d y \, d x
+
\int_{\partial E} \int_{\R^N\setminus E}
K_\delta(x-y)
h(x)^2 H(x) .
\end{align*}

In a similar way, we have
\begin{align*}
A_{22} + B_{22} + C_2
&=
\int_{\partial E} \int_{\partial E}
K_\delta(x-y) h(y)^2 \nu(x) \nu(y)
\, d y \, d x
\\
& \qquad
-
\int_{ E} \int_{\partial E}
K_\delta(x-y)
\Big[ D_{y_i} u^i(y) u^j(y) \nu^j(y)
- u^i(y) D_{y_i}  u^j(y) \nu^j(y) \Big]
\, d y \, d x
\\
&=
\int_{\partial E} \int_{\partial E}
K_\delta(x-y) h(y)^2 \nu(x) \nu(y)
\, d y \, d x
-
\int_{ E} \int_{\partial E}
K_\delta(x-y)
h(y)^2 H(y) \, d y \, d x .
\end{align*}
Further calculations show that
\begin{align*}
A_{12}
&=
-\int_{\partial E} \int_{\partial E}
K_\delta(x-y) h(x) h(y) \, d y d x
\\
&
\qquad
- \int_{\partial E} \int_{\R^N\setminus E}
K_\delta(x-y)  {\rm div}(u)(y) u^i(x) \nu^i(x)
\, d y \, d x
\\
& \qquad
+
\int_{ E} \int_{\partial E}
K_\delta(x-y)  {\rm div}(u)(x) u^i(y) \nu^i(y)
\, d y \, d x
\\
& \qquad
+
\int_{E} \int_{\R^N\setminus E}
K_\delta(x-y)  {\rm div}(u)(x)  {\rm div}(u)(y)
\, d y \, d x ,
\end{align*}
\begin{align*}
A_{21}&=
-\int_{\partial E} \int_{\partial E}
K_\delta(x-y) h(x) h(y) \, d y d x
\\
& \qquad -
\int_{\partial E}\int_{\R^N\setminus E}
K_\delta(x-y)  {\rm div}(u)(y) u^j(x) \nu^j(x)
\, d y \, d x
\\
&\qquad
+
\int_{E} \int_{\partial E}
K_\delta(x-y)  {\rm div}(u)(x) u^i(y) \nu^i(y)
\, d y \, d x
\\
& \qquad
+\int_{E}\int_{\R^N\setminus E}
K_\delta(x-y)  {\rm div}(u)(x) div (u)(y)
\, d y \, d x ,
\end{align*}
and
\begin{align*}
B_{12} + B_{21} &=
2\int_{\partial E} \int_{\R^N\setminus E}
K_\delta(x-y)   {\rm div}(u)(y) u^j(x) \nu^j(x)
\, d y \, d x
\\
& \qquad
-
2\int_{E} \int_{\partial E}
K_\delta(x-y)  {\rm div}(u)(x) u^j(y) \nu^j(y)
\, d y \, d x
\\
& \qquad
-4 \int_{E}\int_{\R^N\setminus E}
K_\delta(x-y)  {\rm div}(u)(x)  {\rm div} (u)(y)
\, d y \, d x ,
\end{align*}
so that
\begin{align*}
A_{12}+A_{21} + B_{12}+B_{21} + C_3
=
-2
\int_{\partial E} \int_{\partial E}
K_\delta(x-y) h(x) h(y) \, d y d x .
\end{align*}
Therefore
\begin{align*}
\frac{d^2}{dt^2}
Per_{s,\delta}(E_{th},\Omega)
\Big|_{t=0}
&=
2\int_{\partial E} \int_{\partial E}
K_\delta(x-y)  h(x)^2  ( \nu(x) \nu(y)-1)
\, d y \, d x
\\
& \qquad
-2
\int_{\partial E}
h(x)
\int_{\partial E}
K_\delta(x-y)  ( h(y) - h(x)) \, d y d x
\\
& \qquad
-
\int_{\partial E}
h(x)^2 H(x)
\int_{\R^N}
(\chi_E(y)- \chi_{E^c}(y))
K_\delta(x-y)
\, d y \, d x .
\end{align*}
Taking the limit as $\delta\to 0$ we find \eqref{second var}.
\qed

\medskip
\noindent
{\bf Proof of proposition~\ref{prop derivative H}.}
Let $\nu_t(x)$ denote the unit normal vector to $\partial E_t$ at $x \in \partial E_t$ pointing out of $E_t$.
Note that $\nu(x) = \nu_0(x)$.
Let $L_t$ be the half space defined by $L_t = \{ x: \langle x - p_t , \nu_t(p_t)\rangle>0\}$.
Then
\begin{align}
\label{int ht}
H_{\Sigma_{th}}^s(p_t) =
\int_{\R^N} \frac{\chi_{E_t} (x)-\chi_{L_t}(x)-\chi_{E^c}(x) + \chi_{L_t^c}(x)} {|x-p_t|^{N+s}} \, d x
\end{align}
since the function $1-2 \chi_{L_t}$ has zero principal value. Note that the integral in \eqref{int ht} is well defined and
\begin{align*}
H_{\Sigma_{th}}^s(p_t)
& =
2
\int_{\R^N} \frac{\chi_{E_t} (x)-\chi_{L_t}(x)} {|x-p_t|^{N+s}} \, d x .
\end{align*}

For $\delta >0$ let $\eta \in C^\infty(\R^N)$ be a radially symmetric cut-off function with $\eta(x)=1 $ for $|x|\geq 2$, $\eta(x) = 0$ for $|x|\leq 1$. Define $\eta_\delta(x) = \eta(x/\delta)$ and
write
$$
\int_{\R^N }
\frac{\chi_{E_t} (x)-\chi_{L_t}(x)} {|x-p_t|^{N+s}} \, d x =
f_\delta (t) + g_\delta(t)
$$
where
\begin{align*}
f_\delta(t)
&=
\int_{\R^N }
\frac{\chi_{E_t} (x)-\chi_{L_t}(x)} {|x-p_t|^{N+s}} \eta_\delta(x-p_t) \, d x
\end{align*}
and $g_\delta(t)$ is the rest.
Then it is direct that $f_\delta$ is differentiable and
\begin{align*}
f_\delta'(0)
&=
\int_{\partial E}
\frac{h(x)}{|x-p|^{N+s}} \eta_\delta(x-p)
\\
& \qquad
-
\int_{\partial L_0}
\frac{h(p) \langle \nu(p) , \nu(p) \rangle
-\langle x-p,\frac{\partial \nu_t(p_t)}{\partial t}|_{t=0}\rangle
}{|x-p|^{N+s}}\eta_\delta(x-p)
\\
&
\qquad
+ (N+s)
h(p)
\int_{\R^N}
\frac{\chi_{E}(x) - \chi_{L_0}(x)}{|x-p|^{N+s+2} }
\langle x-p,\nu(p)\rangle \eta_\delta(x-p) dx
\\
&
\qquad
-
h(p)
\int_{\R^N}
\frac{\chi_{E}(x) - \chi_{L_0}(x)}{|x-p|^{N+s} }
\langle \nabla \eta_\delta(x-p),\nu(p)\rangle  dx .
\end{align*}
We integrate the third term by parts
\begin{align*}
& (N+s)
\int_{\R^N}
\frac{\chi_{E}(x) - \chi_{L_0}(x)}{|x-p|^{N+s+2} }
\langle x-p,\nu(p)\rangle \eta_\delta(x-p) dx
\\
&=
-
\int_{\R^N}
( \chi_{E}(x) - \chi_{L_0}(x) )
\langle \nabla\frac{1}{|x-p|^{N+s} }
,\nu(p)\rangle \eta_\delta(x-p) dx
\\
& =
-
\int_{\partial E} \frac{\langle \nu(x) ,\nu(p) \rangle }{|x-p|^{N+s}}
\eta_\delta(x-p)
+
\int_{\partial L_0}
\frac{\langle \nu(p), \nu(p) \rangle }{|x-p|^{N+s}}
\eta_\delta(x-p)
\\
& \qquad +\int_{\R^N}\frac{\chi_{E(x)} - \chi_{L_0}(x)}{|x-p|^{N+s}}
\langle \nabla \eta_\delta(x-p) , \nu(p) \rangle dx .
\end{align*}
Since $\eta_\delta $ is radially symmetric,
$$
\int_{\partial L_0}
\frac{
\langle x-p, \frac{\partial \nu_t(p_t)}{\partial t}|_{t=0} \rangle
}{|x-p|^{N+s}}\eta_\delta(x-p) \, d x=0
$$
and then
\begin{align*}
f_\delta'(0)
&=
\int_{\partial E}
\frac{h(x)}{|x-p|^{N+s}} \eta_\delta(x-p)
dx
-
h(p)
\int_{\partial E} \frac{\langle \nu(x) ,\nu(p) \rangle }{|x-p|^{N+s}}
\eta_\delta(x-p) dx,
\end{align*}
which we write as
\begin{align*}
f_\delta'(0)
&=
\int_{\partial E}
\frac{h(x)-h(p)}{|x - p|^{N+s}} \eta_\delta(x-p)
\, d x
+
h(p)
\int_{\partial E} \frac{1-\langle \nu(x) ,\nu(p) \rangle }{|x-p|^{N+s}}
\eta_\delta(x-p) \, d x.
\end{align*}
We claim that  $g_\delta'(t)\to 0$ as $\delta \to 0$, uniformly for $t$ in a neighborhood of 0.
Indeed, in a neighborhood of $p_t$ we can represent $\partial E_t$ as a graph of a function $G_t$ over $L_t \cap B(p_t,2\delta)$, with $G_t$ defined in a neighborhood of $0$ in $\R^{N-1}$, $G_t(0)=0$, $\nabla_{y'} G_t(0)=0$ and smooth in all its variables (we write $y' \in \R^{N-1}$).
Then $g_\delta(t)$ becomes
$$
g_\delta(t)
=
\int_{|y'|<2\delta}
\int_0^{G_t(y')}
\frac{1}{(|y'|^2 + y_N^2)^{\frac{N+s}{2}}}
(1-\eta_\delta(y', y_N) )
d y_N
d y'
$$
so that
$$
g_\delta'(t)
=
\int_{|y'|<2\delta}
\frac{1}{(|y'|^2 + G_t(y')^2)^{\frac{N+s}{2}}}
\frac{\partial G_t}{\partial t}(y')
(1-\eta_\delta(y', y_N) )
d y' .
$$
But $|G_t(y')|\leq K |y'|^2$ and $| \frac{\partial G_t}{\partial t} (y') |\leq K |y'|^2$, so
$$
g_\delta'(t) \leq C \delta^{1-s} .
$$
Therefore
\begin{align*}
\frac{d}{d t}
H_{\Sigma_{th}}^s(p_t)
\Big|_{t=0}
& = 2
\lim_{\delta \to 0}
\Bigg[
\int_{\partial E}
\frac{h(x)-h(p)}{|x - p|^{N+s}} \eta_\delta(x-p) dx
\\
& \qquad
+
h(p)
\int_{\partial E} \frac{1-\langle \nu(x) ,\nu(p) \rangle }{|x-p|^{N+s}}
\eta_\delta(x-p) dx
\Bigg] .
\end{align*}
Letting $\delta\to 0$ we find \eqref{der mean C}.
\qed

\medskip
\noindent
{\bf Proof of corollary~\ref{coro entire stable}.}
The same argument as in the proof of proposition~\ref{prop derivative H} shows that if $F:\Sigma \to \R^N$ is a smooth bounded vector field and we let
$E_{t}$ be the set whose boundary $\Sigma_t = \pp E_{t} $ is parametrized as
$$
\pp E_{th} =  \{ x+ t F(x)\ /\ x\in \pp E \},
$$
with exterior normal vector close to $\nu$, then
$$
\frac{d}{d t}
H_{\Sigma_t}^s(p_t)
\Big|_{t=0}
= 2 \JJ_\Sigma^s[ \langle F , \nu\rangle] (p) ,
$$
where $p_t = p + t F(p)$.
Taking as $F(x) = e_N = (0,\ldots,0,1)$ we conclude that $w = \langle \nu,e_N \rangle$ is a positive function satisfying
$$
\JJ_\Sigma^s[ w ] (x) =0
 \quad\text{for all } x \in \Sigma.
$$
More explicitly
\begin{align}
\label{20a}
\text{p.v.}
\int_{\Sigma}
\frac{w(y)-w(x)}{|y-x|^{N+s}} dy
+
w(x) A(x) = 0 \quad\text{for all } x \in \Sigma,
\end{align}
where
$$
A(x) = \int_\Sigma
\frac{\langle \nu (x)- \nu(y), \nu (x) \rangle }{|x-y|^{N+s}} dy.
$$

As in the classical setting we can show that $\Sigma$ is stable in the sense that \eqref{linearly stable} holds.
Let $\phi \in C_0^\infty(\Sigma)$ and
observe that
\begin{align*}
\frac12
\int_{\Sigma}\int_{\Sigma}
\frac{(\phi(x) - \phi(y))^2}{|x-y|^{N+s} } d x d y
=
\int_{\Sigma}\int_{\Sigma}
\frac{(\phi(x) - \phi(y) ) \phi(x)}{|x-y|^{N+s} } d x d y .
\end{align*}
Write $\phi = w \psi$ with $\psi \in C_0^\infty(\Sigma)$. Then
\begin{align}
\nonumber
\int_{\Sigma}\int_{\Sigma}
\frac{(\phi(x) - \phi(y) ) \phi(x)}{|x-y|^{N+s} } d x d y
& =
\int_{\Sigma}\int_{\Sigma}
\frac{(w(x) - w(y) ) w(x) \psi(x)^2}{|x-y|^{N+s} } d x d y
\\
\label{term2}
& \quad +
\int_{\Sigma}\int_{\Sigma}
\frac{(\psi(x) - \psi(y) ) w(x) w(y) \psi(x)}{|x-y|^{N+s} } d x d y .
\end{align}
Multiplying \eqref{20a} by $w \psi^2$ and integrating we get
\begin{align}
\label{term3}
\int_{\Sigma}\int_{\Sigma}
\frac{(w(x) - w(y) ) w(x) \psi(x)^2}{|x-y|^{N+s} } d x d y
=
\int_{\Sigma}A(x) w(x)^2 \psi(x)^2 dx
=
\int_{\Sigma} A(x) \phi(x)^2 dx .
\end{align}
For the second term in \eqref{term2} we observe that
\begin{align}
\label{term4}
\int_{\Sigma}\int_{\Sigma}
\frac{(\psi(x) - \psi(y) ) w(x) w(y) \psi(x)}{|x-y|^{N+s} } d x d y
 & =
\frac12
\int_{\Sigma}\int_{\Sigma}
\frac{(\psi(x) - \psi(y) )^2 w(x) w(y) }{|x-y|^{N+s} } d x d y .
\end{align}
Therefore, combining \eqref{term2}, \eqref{term3}, \eqref{term4} we obtain
\begin{align*}
\frac12
\int_{\Sigma}\int_{\Sigma}
\frac{(\phi(x) - \phi(y))^2}{|x-y|^{N+s} } d x d y
& =
\int_{\Sigma}A(x)\phi(x)^2 dx
\\
\nonumber
& \quad +
\frac12
\int_{\Sigma}\int_{\Sigma}
\frac{(\psi(x) - \psi(y) )^2 w(x) w(y) }{|x-y|^{N+s} } d x d y .
\end{align*}
and tis shows \eqref{linearly stable}.
\qed

\section{Stability and instability}
\label{sect stability}

We consider the nonlocal minimal cone $C_m^n(s) =  \partial E_\alpha$
where $E_\alpha$ is defined in \eqref{def solid cone} and $\alpha$ is the one of theorem~\ref{thm1}.
For $0 \leq s <1$ we obtain a characterization of
their stability in terms of constants that depend on   $m$, $n$ and $s$. For the case $s=0$ we consider the limiting cone with parameter $\alpha_0$ given in proposition~\ref{prop alpha s=0}.
Note that in the case $s=0$ the limiting Jacobi operator $\mathcal J^0_{C_{\alpha_0}}$ is well defined for smooth functions with compact support.

For brevity, in this section we write $\Sigma = C_m^n(s)$.

Recall that
$$
\mathcal J^s_{\Sigma}[\phi](x)
=
\text{p.v.}
\int_{\Sigma}
\frac{\phi(y) - \phi(x)}{|y-x|^{N+s}} dy
+ \phi(x) \int_{\Sigma}
\frac{1-\langle \nu(x),\nu(y)\rangle}{|x-y|^{N+s}} dy
$$
for $\phi \in C_0^\infty(\Sigma\setminus\{0\})$.
Let us rewrite this operator in the form
$$
\mathcal J^s_{\Sigma}[\phi](x)
=
\text{p.v.}
\int_{\Sigma}
\frac{\phi(y) - \phi(x)}{|x-y|^{N+s}} dy
+ \frac{A_0(m,n,s)^2}{|x|^{1+s}}
\phi(x)
$$
where
\begin{align*}
A_0(m,n,s)^2 = \int_{\Sigma}
\frac{\langle \nu (\hat p)- \nu(x), \nu (\hat p) \rangle }{|\hat p-x|^{N+s}} dx \geq0
\end{align*}
and this integral is evaluated at any $\hat p\in {\Sigma}$ with $|\hat p|=1$.
We can think of $\mathcal J^s_{\Sigma}$ as analogous to the fractional Hardy operator
$$
-(-\Delta)^{\frac{1+s}2} \phi
+ \frac c{|x|^{1+s}} \phi
\quad \text{in } \R^{N-1} ,
$$
for which positivity is related to a fractional Hardy inequality
with best constant, see Herbst \cite{herbst}.
This suggests that the positivity of $\mathcal J_\Sigma$ is related to the existence of $\beta$ in an appropriate range such that
$\mathcal J^s_{\Sigma}[|x|^{-\beta}]\leq 0$, and it turns out that the best choice of $\beta$ is $\beta = \frac{N-2-s}{2}$.
This motivates the definition
$$
H(m,n,s) =
\text{p.v.}
\int_{\Sigma}
\frac{1 - |y|^{-\frac{N-2-s}{2}}}{|\hat p-y|^{N+s}} dy
$$
where $\hat p\in \Sigma$ is any point with $|\hat p|=1$.



We have then the following Hardy inequality with best constant:
\begin{prop}
\label{prop hardy ineq}
For any $\phi \in C_0^\infty(\Sigma\setminus\{0\})$ we have
\begin{align}
\label{fract hardy}
H(m,n,s)
\int_{\Sigma}
\frac{\phi(x)^2}{|x|^{1+s}} dx
\leq
\frac{1}{2}
\int_{\Sigma}
\int_{\Sigma}
\frac{(\phi(x)-\phi(y) ) ^2}{|x-y|^{N+s}} d x dy
\end{align}
and $H(m,n,s)$ is the best possible constant in this inequality.
\end{prop}

As a result we have:
\begin{corollary} The cone $C_m^n(s)$ is stable if and only if
$H (m,n,s) \geq A_0(m,n,s)^2$.
\end{corollary}

Other related fractional Hardy inequalities have appeared in the literature, see for instance \cite{bogdan-dyda,dyda-frank}.

\medskip
\noindent
{\bf Proof of proposition~\ref{prop hardy ineq}.}
Let us write  $H = H(m,n,s)$ for simplicity.
To prove the validity of \eqref{fract hardy}
let $w(x) = |x|^{-\beta}$ with $\beta = \frac{N-2-s}{2}$ so that from the definition of $H$ and homogeneity we have
\begin{align*}
\text{p.v.}
\int_{\Sigma} \frac{w(y) - w(x)}{|y-x|^{N+s} } d y
+ \frac{H}{|x|^{1+s}} w(x) = 0
\quad \text{for all } x \in\Sigma\setminus \{0\}.
\end{align*}
Now the same argument as in the proof of corollary~\ref{coro entire stable} shows that
\begin{align}
\label{equiv form}
\frac12
\int_{\Sigma}\int_{\Sigma}
\frac{(\phi(x) - \phi(y))^2}{|x-y|^{N+s} } d x d y
& =
\int_{\Sigma}\frac{H}{|x|^{1+s}}\phi(x)^2 dx
\\
\nonumber
& \quad +
\frac12
\int_{\Sigma}\int_{\Sigma}
\frac{(\psi(x) - \psi(y) )^2 w(x) w(y) }{|x-y|^{N+s} } d x d y .
\end{align}
for all $\phi \in C_0^\infty(\Sigma\setminus\{0\})$  with $\psi = \frac\phi w \in C_0^\infty(\Sigma\setminus\{0\})$

Now let us show that $H$ is the best possible constant in \eqref{fract hardy}. Assume that
$$
\tilde H
\int_{\Sigma}
\frac{\phi(x)^2}{|x|^{1+s}} dx
\leq
\frac{1}{2}
\int_{\Sigma}
\int_{\Sigma}
\frac{(\phi(x)-\phi(y) ) ^2}{|x-y|^{N+s}} d x dy
$$ for all
$\phi \in C_0^\infty(\Sigma \setminus\{0\})$.
Using \eqref{equiv form} and letting $\phi = w \psi$ with $\psi\in \in C_0^\infty(\Sigma \setminus\{0\})$ we then have
\begin{align*}
\tilde H
\int_{\Sigma}
\frac{w(x)^2 \psi(x)^2}{|x|^{1+s}} dx
& \leq
H \int_{\Sigma}
\frac{w(x)^2 \psi(x)^2}{|x|^{1+s}} dx
\\
& \qquad +\frac12
\int_{\Sigma}\int_{\Sigma}
\frac{(\psi(x) - \psi(y) )^2 w(x) w(y) }{|x-y|^{N+s} } d x d y .
\end{align*}
For $R>3$ let $\psi_R:\Sigma\to [0,1]$ be a radial function such that $\psi_R(x)=0$ for $|x|\leq 1$, $\psi_R(x)=1$ for $2\leq |x|\leq 2 R$, $\psi_R(x)=0$ for $|x|\geq 3 R$. We also require $|\nabla\psi_R(x)|\leq C$ for $|x|\leq 3$, $|\nabla\psi_R(x)|\leq C/R$ for $2 R\leq |x|\leq 3 R$. We claim that
\begin{align}
\label{ineq 1}
a_0 \log(R) - C
\leq
\int_{\Sigma}
\frac{w(x)^2 \psi_R(x)^2}{|x|^{1+s}} dx \leq a_0 \log(R) +C
\end{align}
where $a_0>0$, $C>0$ are independent of $R$, while
\begin{align}
\label{ineq 2}
\left|
\int_{\Sigma}\int_{\Sigma}
\frac{(\psi_R(x) - \psi_R(y) )^2 w(x) w(y) }{|x-y|^{N+s} } d x d y
\right|\leq C .
\end{align}
Letting then $R\to \infty$ we deduce that $\tilde H \leq H$.

To prove the upper bound in \eqref{ineq 1} let us write points in $\Sigma$ as $x = (y,z) $, with $y \in \R^m$, $z\in\R^n$.
Let us write
$ y = r \omega_1 $, $z = r\omega_2$, with $r>0$,  $\omega_1 \in S^{m-1}$, $\omega_2 \in S^{n-1}$ and use spherical coordinates $(\theta_1,\ldots,\theta_{m-1})$ and $(\varphi_1,\ldots,\varphi_{n-1})$ for $\omega_1$ and $\omega_2$ as in \eqref{omega1} and \eqref{omega2} . We assume here that $m\geq n \geq 2 $. In the remaining cases the computations are similar. Then we have
$$
\int_{\Sigma}
\frac{w(x)^2 \psi_R(x)^2}{|x|^{1+s}} dx \leq
a_0
\int_1^{4 R} \frac{1}{r^{N-2-s}} \frac1{r^{1+s}} r^{N-2} dr
\leq a_0 \log(R) + C
$$
where
$$
a_0 = \sqrt{1+\alpha^2} A_{m-1} A_{n-1}
$$
and $A_k$ denotes the area of the sphere $S^k \subseteq \R^{k+1}$ and is given by
\begin{align}
\label{Am}
A_{k} = \frac{2\pi^{\frac{k+1}2}}{\Gamma(\frac{k+1}2)}.
\end{align}
The lower bound in \eqref{ineq 1} is similar.

To obtain \eqref{ineq 2} we split $\Sigma$ into the regions $R_1 = \{x: |x|\leq 3\}$, $R_2 = \{x : 3\leq x \leq R\}$, $R_3=\{x: R\leq |x|\leq4R\}$ and $R_4 = \{x:|x|\geq 4R\}$ and let
$$
I_{i,j} = \int_{x\in R_i} \int_{y\in R_j}
\frac{(\psi_R(x) - \psi_R(y) )^2 w(x) w(y) }{|x-y|^{N+s} } d x d y .
$$
Then $I_{i,j} = I_{j,i}$ and $I_{j,j}=0$ for $j=2,4$. Moreoover $I_{1,1}=O(1)$ since the region of integration is bounded and $\psi_R$ is uniformly Lipschitz.

Estimate of $I_{1,2}$: We bound $w(x)\leq C$ for $|x|\geq 1$ and then
\begin{align*}
|I_{1,2}|\leq C  \int_{y\in R_2} \frac{w(y)}{|p-y|^{N+s}} dy
\leq C \int_2^R \frac{1}{r^{\frac{N-2-s}{2}}} \frac1{r^{N+s}} r^{N-2} dr \leq C,
\end{align*}
where  $p\in\Sigma$ is fixed with $|p|=2$.

By the same argument  $I_{1,3}=O(1) $ and $I_{1,4}=O(1)$ as $R\to\infty$.

Estimate of $I_{2,3}$: for $y \in R_3$, $w(y)\leq C R^{-\frac{N-2-s}2}$, so
\begin{align*}
|I_{2,3}|
&\leq
C R^{-\frac{N-2-s}2}
\int_{x\in R_2} \frac{1}{|x|^{\frac{N-2-s}2}}
\int_{y\in R_3}
\frac{(\psi_R(x)-\psi_R(y))^2}{|x-y|^{N+s}} dy dx
\\
&
\leq
C R^{-\frac{N-2-s}2}
\frac{Vol(R_3)}{R^{N+s}} \int_{x\in R_2}
\frac{1}{|x|^{\frac{N-2-s}2}}  dx \leq C.
\end{align*}

Estimate of $I_{2,4}$:
\begin{align*}
|I_{2,4}|
\leq C
\int_{x\in R_2}
\frac{1}{|x|^{\frac{N-2-s}2}}
\int_{y\in R_4}
\frac{1}{|x-y|^{N+s}}
\frac{1}{|y|^{\frac{N-2-s}{2}}}
dy
dx .
\end{align*}
By scaling
\begin{align*}
\int_{y\in R_4}
\frac{1}{|x-y|^{N+s}}
\frac{1}{|y|^{\frac{N-2-s}{2}}}
dy \leq C R^{-\frac N2  - \frac s2}
\quad \text{for } x \in R_2 ,
\end{align*}
so that
\begin{align*}
|I_{2,4}|
\leq C R^{-\frac N2  - \frac s2}
\int_{x\in R_2}
\frac{1}{|x|^{\frac{N-2-s}2}}
dx \leq C  .
\end{align*}

To estimate $I_{3,3}$ we use $|\psi_R(x) -\psi_R(y)|\leq \frac CR |x-y$ for $x,y \in R_3$, which yields
\begin{align*}
|I_{3,3}|
&\leq \frac{C}{R^2} \frac{1}{R^{N-2-s}} \int_{x,y\in R_3} \frac{1}{|x-y|^{N+s-2}}
dy dx .
\end{align*}
The integral is finite and by scaling we see that is bounded by $C R^{N-s}$, so that
$$
|I_{3,3}|\leq C.
$$

Estimate of $I_{3,4}$:
\begin{align*}
|I_{3,4}|
\leq C R^{-\frac{N-2-s}2}
\int_{x\in R_3}
\int_{y\in R_4}
\frac{1}{|x-y|^{N+s}} \frac{1}{|y|^{\frac{N-2-s}{2}}} d y d x .
\end{align*}
By scaling
$$
\int_{y\in R_4}
\frac{1}{|x-y|^{N+s}} \frac{1}{|y|^{\frac{N-2-s}{2}}} d y
\leq \frac{C}{|x|^{\frac{N+s}2}}
$$
for $x\in R_3$. Therefore
$$
|I_{3,4}|
\leq C R^{-\frac{N-2-s}2}
\int_{x\in R_3} \frac{1}{|x|^{\frac{N+s}2}} dx
\leq C.
$$
This concludes the proof of \eqref{ineq 2}.
\qed

\bigskip

\noindent
{\bf Proof of Theorem~\ref{thm stability}.}
In what follows we will obtain expressions for $H(m,n,s)$ and $A_0(m,n,s)^2$ for $m\geq 2$, $n\geq 1$, $0\leq s <1$. We always assume $m\geq n$.
For the sake of generality, we will compute
$$
C(m,n,s,\beta) =
\text{p.v.}
\int_{\Sigma}
\frac{1 - |x|^{-\beta}}{|\hat p-x|^{N+s}} dx
$$
where $\hat p\in \Sigma$, $|\hat p|=1$, and $\beta \in (0,N-2-s)$, so that $H(m,n,s) = C(m,n,s,\frac{N-2-s}{2})$.

Let $x = (y,z) \in \Sigma$, with $y \in \R^m$, $z\in\R^n$.
For simplicity in the next formulas we take $p = (e_1^{(m)},\alpha e_2^{(n)})$ (see the notation in \eqref{notation e}), and $h(y,z) = |y|^{-\beta}$,
so that
$$
C(m,n,s,\beta)
=
(1+\alpha^2)^{\frac{1+s}2}
\text{p.v.}
\int_{\Sigma}
\frac{h(p) - h(x)}{|p-x|^{N+s}} \, d x .
$$

\medskip

\noindent
{\bf Computation of $C(m,1,s,\beta)$.}
Write
$ y = r \omega_1 $, $z = \pm \alpha r$, with $r>0$,  $\omega_1 \in S^{m-1}$.
Let us use the notation $\Sigma_{\alpha}^+ = \Sigma \cap [ z>0]$, $\Sigma_{\alpha}^- = \Sigma \cap [ z<0]$.
Using polar coordinates $(\theta_1,\ldots,\theta_{m-1})$
for $\omega_1$  as in \eqref{omega1}  we have
$$
|x-p|^2 = |r\theta_1-e_1^{(m)}|^2 + \alpha^2 |r\theta_1-e_1^{(m)}|^2
=  r^2 + 1 - 2 r \cos(\theta_1)  + \alpha^2 (r-1)^2,
$$
for $x\in \Sigma_{\alpha}^+$ and
$$
|x-p|^2 = |r\theta_1-e_1^{(m)}|^2 + \alpha^2 |r\theta_1-e_1^{(m)}|^2
=  r^2 + 1 - 2 r \cos(\theta_1)  + \alpha^2 (r+1)^2,
$$
for $x\in \Sigma_{\alpha}^-$.
Hence, with $h(y,z) = |y|^{-\beta}$
\begin{align}
\label{int n1 h}
\text{p.v.}
\int_{\Sigma}
\frac{h(p)- h(x)}{|x-p|^{N+s}} dx
=
\sqrt{1+\alpha^2}
A_{m-2}
\text{p.v.}
\int_0^\infty (1-r^{-\beta})
( I_+(r) + I_-(r) )
r^{N-2}
dr
\end{align}
where
\begin{align*}
I_+(r) & =
\int_0^\pi
\frac{\sin(\theta_1)^{m-2} }
{(r^2 + 1 - 2 r \cos(\theta_1)  + \alpha^2 (r-1)^2)^{\frac{N+s}2}}
d\theta_1
\\
I_-(r) & =
\frac{\sin(\theta_1)^{m-2} }
{(r^2 + 1 - 2 r \cos(\theta_1)  + \alpha^2 (r+1)^2)^{\frac{N+s}2}}
d\theta_1 ,
\end{align*}
and $A_{m-2}$ is defined in \eqref{Am} for $m\geq 2$.
From \eqref{int n1 h}
we obtain
\begin{align}
\label{Cm1sbeta}
C(m,1,s,\beta)
=
(1+\alpha^2)^{\frac{3+s}2}
A_{m-2}
\int_0^1 (r^{N-2}-r^{N-2-\beta}+r^s- r^{\beta+s})
( I_+(r) + I_-(r) )
d r .
\end{align}

\medskip
\noindent
{\bf Computation of $A_0(m,1,s)^2$.}
Let $x = (r \theta_1,\pm \alpha r) $, $p = (e_1^{(n)},\alpha)$
so that
\begin{align*}
\nu(x ) = \frac{ ( -\alpha  \omega_1, \pm 1)}{\sqrt{1+\alpha^2}} ,
\quad
\nu(p) = \frac{ ( -\alpha e_1^{(n)}, 1)}{\sqrt{1+\alpha^2}} ,
\end{align*}
and hence
\begin{align*}
\int_{\Sigma}
\frac{1-\langle \nu(x), \nu (p) \rangle }{|p-x|^{N+s}} dx
&=
\sqrt{1+\alpha^2}
A_{m-2}
\int_0^\infty
( J_+(r) + J_-(r)) r^{N-2} dr
\\
&=
\sqrt{1+\alpha^2}
A_{m-2}
\int_0^1
(r^{N-2} + r^s)
( J_+(r) + J_-(r))  dr ,
\end{align*}
where
\begin{align*}
J_+(r)
&=
\frac{\alpha^2}{1+\alpha^2}
\int_0^\pi
\frac{(1-\cos(\theta_1)) \sin(\theta_1)^{m-2} }
{(r^2 + 1 - 2 \cos(\theta_1)  + \alpha^2 (r-1)^2)^{\frac{N+s}2}}
d\theta_1
\\
J_-(r)
&=
\frac{1}{1+\alpha^2}
\int_0^\pi
\frac{[ 2 +\alpha^2-\alpha^2 \cos(\theta_1))\sin(\theta_1)^{m-2} }
{(r^2 + 1 - 2 r \cos(\theta_1)  + \alpha^2 (r+1)^2)^{\frac{N+s}2}}
d\theta_1
\end{align*}
Therefore we find
\begin{align*}
A_0(m,1,s)^2
& =
(1+\alpha^2)^{\frac{3+s}2}
A_{m-2}
\int_0^1
(r^{N-2} + r^s)
( J_+(r) + J_-(r))  dr .
\end{align*}

\noindent
{\bf Computation of $C(m,n,s,\beta)$ for $n\geq 2$.}
Write
$ y = r \omega_1 $, $z = r\omega_2$, with $r>0$,  $\omega_1 \in S^{m-1}$, $\omega_2 \in S^{n-1}$ and let us use spherical coordinates $(\theta_1,\ldots,\theta_{m-1})$ and $(\varphi_1,\ldots,\varphi_{n-1})$ for $\omega_1$ and $\omega_2$ as in \eqref{omega1}  and \eqref{omega2}. Recalling that
$p = (e_1^{(m)},\alpha e_2^{(n)})$, we have
$$
|x-p|^2 = |r\theta_1-e_1^{(m)}|^2 +|r\theta_1-e_1^{(m)}|^2
=  r^2 + 1 - 2 r \cos(\theta_1)  + \alpha^2 ( r^2 + 1 - 2 r \cos(\varphi_1)).
$$
Hence, with $h(y,z) = |y|^{-\beta}$
\begin{align*}
\text{p.v.}
&\int_{\Sigma}
\frac{h(p)- h(x)}{|x-p|^{N}} dx
=
\sqrt{1+\alpha^2}
A_{m-2}A_{n-2}
\text{p.v.}
\int_0^\infty (1-r^{-\beta})
I(r) r^{N-2}  dr
\\
&=
\sqrt{1+\alpha^2}
A_{m-2}A_{n-2}
\int_0^1 (r^{N-2}-r^{N-2-\beta}+r^s- r^{\beta+s})
I(r)
d r
\end{align*}
where
$$
I(r)
=
\int_0^\pi
\int_0^\pi \frac{\sin(\theta_1)^{m-2} \sin(\varphi_1)^{n-2}}{(r^2 + 1 - 2 r \cos(\theta_1)  + \alpha^2 ( r^2 + 1 - 2 r \cos(\varphi_1)))^{\frac{N+s}{2}}}
d\theta_1 d\varphi_1.
$$
We find then that
\begin{align}
\label{Cmnsbeta}
C(m,n,s,\beta) =
(1+\alpha)^{\frac{3+s}{2}}
A_{m-2}A_{n-2}
\int_0^1 (r^{N-2}-r^{N-2-\beta}+r^s- r^{\beta+s})
I(r)
d r .
\end{align}
\medskip
\noindent
{\bf Computation of $A_0(m,n,s)^2$ for $n\geq 2$.}
Similarly as before we have, for $x=(r \omega_1,\alpha r \omega_2)\in \Sigma$, and $p=(e_1^{(m)},\alpha e_2^{(n)} )$:
\begin{align*}
\nu(x ) = \frac{ ( -\alpha \omega_1, \omega_2)}{\sqrt{1+\alpha^2}} ,
\quad
\nu(p) = \frac{ ( -\alpha e_1^{(n)}, 1)}{\sqrt{1+\alpha^2}} .
\end{align*}
Hence
\begin{align*}
\int_{\Sigma}
\frac{1-\langle \nu(x), \nu (p) \rangle }{|p-x|^{N+s}} dx
& =
\sqrt{1+\alpha^2}
A_{m-2}
A_{n-2}
\int_0^\infty r^{N-2} J(r) dr
\\
&=
\sqrt{1+\alpha^2}
A_{m-2}
A_{n-2}
\int_0^1(r^{N-2} + r^s) J(r) dr
\end{align*}
where
\begin{align*}
J(r) = \frac{1}{1+\alpha^2}
\int_0^\pi
\int_0^\pi
\frac{(1+\alpha^2 -\alpha^2 \cos(\theta_1) - \cos(\varphi_1) ) \sin(\theta_1)^{m-2} \sin(\varphi_1)^{n-2}}{ ( r^2 + 1 - 2r\cos(\theta_1) + \alpha^2 ( r^2 + 1 - 2 r \cos(\varphi_1)
)^{\frac{N+s}{2}}} d\theta_1 d\varphi_1 .
\end{align*}
We finally obtain
\begin{align*}
A_0(m,n,s)^2
= (1+\alpha^2)^{\frac{3+s}{2}}
A_{m-2}
A_{n-2}
\int_0^1(r^{N-2} + r^s) J(r) dr.
\end{align*}

In table~1 we show the  values obtained for $H(m,n,0)$ and $A_0(m,n,0)^2$, divided by $(1+\alpha^2)^{\frac{3+s}{2}}
A_{m-2}
A_{n-2}$, from numerical approximation of the integrals.
From these results we can say that for $s=0$, $\Sigma$ is stable if $n+m=7$ and unstable if $n+m\leq 6$.
The same holds for $s>0$ close to zero by continuity of the values with respect to $s$.
\qed

\begin{table}

\begin{tabular}{|c|l|l|l|l|l|l|l|l|}
\hline
&&\multicolumn{6}{c}{$n$}&
\\
\hline
 &  & 1 & 2 & 3 & 4 & 5 & 6 & 7
\\
$m$ & & & & & & & &
\\
\hline
2 & $H $  & 0.8140 & 1.0679 & & & & &
\\
 & $A_0^2$  & 3.2669 & 2.3015 & & & & &
\\
\hline
3 & $H $ & 1.1978 & 1.2346 & 0.3926 & & & &
\\
 & $A_0^2$ & 2.5984 & 1.7918 & 0.4463 & & & &
\\
\hline
4 & $H $ & 1.3968 & 1.3649 & 0.4477 & 0.1613 & & &
\\
 & $A_0^2$ & 2.0413 & 1.5534 & 0.4288 & 0.1356 & & &
\\
\hline
5 & $H $ & 1.5117 & 1.4570 & 0.4895 & 0.1845 & 0.06978 & &
\\
 & $A_0^2$ & 1.7332 & 1.3981 &0.4118 & 0.1398 & 0.04849 & &
\\
\hline
6 & $H $ & 1.5833 & 1.5231 & 0.5215 & 0.2031 & 0.08013 & 0.03113 &
\\
 & $A_0^2$ & 1.5318 & 1.2841 & 0.3955 &  0.1412 & 0.05173 &0.01885 &
\\
\hline
7 & $H $ & 1.6303 & 1.5719 & 0.5465 & 0.2182 & 0.08885 & 0.03583 & 0.01416
\\
& $A_0^2$ & 1.3872 & 1.1951 & 0.3802 & 0.1409 & 0.05381 & 0.02051 &0.007704
\\
\hline

\end{tabular}

\bigskip

\caption{ Values of $H(m,n,0)$ and $A_0(m,n,0)^2$ divided by $(1+\alpha^2)^{\frac{3+s}{2}}
A_{m-2}
A_{n-2}$}
\end{table}

\begin{remark}
We see from formulas \eqref{Cm1sbeta} and \eqref{Cmnsbeta} that $C(m,n,s,\beta)$ is symmetric with respect to $\frac{N-2-s}{2}$ and is maximized for $\beta = \frac{N-2-s}{2}$.
\end{remark}

\begin{remark}
In table 2 we give some numerical values of $\alpha$, $H(m,n,s)$ and $A_0(m,n,s)^2$ divided by $(1+\alpha^2)^{\frac{3+s}{2}}
A_{m-2}
A_{n-2}$ for $m=4$, $n=3$, which show how in this dimension stability depends on $s$. One may conjecture that there is $s_0$ such that the cone is stable for $0 \leq s \leq s_0$ and unstable for $s_0<s<1$.
\end{remark}

\begin{table}

\begin{tabular}{|l|l|l|l|l|}
\hline
&\multicolumn{3}{c}{$s$}&
\\
\hline
& 0.1 & 0.2 & 0.3 & 0.4
\\
\hline
$\alpha$ & 0.8379 & 0.8361 & 0.8341	& 0.8319
\\
$H(4,3,s)$ & 0.4113 & 0.3856 & 0.3699 & 0.3639
\\
$A_0(4,3,s)^2$ & 0.4007 & 0.3830 & 0.3756 & 0.3786
\\
\hline
\end{tabular}

\bigskip

\caption{Values of $H(m,n,s)$ and $A_0(m,n,s)^2$ divided by $(1+\alpha^2)^{\frac{3+s}{2}}
A_{m-2}
A_{n-2}$ for $m=4$, $n=3$.}

\end{table}

\end{document}